\documentclass{article}

\usepackage{graphicx}

\newtheorem{thm}{Theorem}[section]
\newtheorem{cor}[thm]{Corollary}
\newtheorem{lem}[thm]{Lemma}
\newtheorem{de}[thm]{Definition}
\newcounter{bean}
\newcounter{milk}

\begin{document}

\title{The Mystery of the Shape Parameter IV}         
\author{Lin-Tian Luh\\Department of Mathematics, Providence University\\ Shalu Town, Taichung County\\ Taiwan\\Email: ltluh@pu.edu.tw}        
\date{\today}          
\maketitle

{\bf Abstract}. This is the fourth paper of our study of the shape parameter c contained in the famous multiquadrics $(-1)^{\lceil \beta\rceil}(c^{2}+\|x\|^{2})^{\beta},\ \beta>0$, and the inverse multiquadrics $(c^{2}+\|x\|^{2})^{\beta},\ \beta<0$. The theoretical ground is the same as that of \cite{Lu6}. However we extend the space of interpolated functions to a more general one. This leads to a totally different set of criteria of choosing c.\\
\\
{\bf keywords}: radial basis function, multiquadric, shape parameter, interpolation.

\section{Introduction}       
Again, we are going to adopt the radial function
\begin{eqnarray}
h(x):=\Gamma(-\frac{\beta}{2})(c^{2}+|x|^{2})^{\frac{\beta}{2}},\ \beta\in R\backslash 2N_{\geq 0},\ c>0
\end{eqnarray}
, where $|x|$ is the Euclidean norm of $x$ in $R^{n},\ \Gamma$ is the classical gamma function, and $c,\beta$ are constants. This definition looks more complicated than the ones mentioned in the abstract. However it will simplify the Fourier transform of $h$ and our analysis of some useful results.

In order to make this paper more readable, we review some basic ingredients mentioned in the previous papers, at the cost of wasting a few pages.

For any interpolated function $f$, our interpolating function will be of the form
\begin{eqnarray}
s(x):=\sum_{i=1}^{N}c_{i}h(x-x_{i})+p(x)
\end{eqnarray}
where $p(x)\in P_{m-1}$, the space of polynomials of degree less than or equal to $m-1$ in $R^{n}, X=\{x_{1},\cdots,x_{N}\}$ is the set of centers(interpolation points). For $m=0,\ P_{m-1}:=\{0\}$. We require that $s(\cdot )$ interpolate $f(\cdot )$ at data points $(x_{1},f(x_{1})),\cdots,(x_{N},f(x_{N}))$. This results in a linear system of the form
\begin{eqnarray}
  \sum_{i=1}^{N}c_{i}h(x_{j}-x_{i})+\sum_{i=1}^{Q}b_{i}p_{i}(x_{j})=f(x_{j}) &  & ,j=1,\cdots,N \nonumber \\
          \\
  \sum_{i=1}^{N}c_{i}p_{j}(x_{i})=0                                          &  & ,j=1,\cdots,Q \nonumber
\end{eqnarray}
to be solved, where $\{p_{1},\cdots,p_{Q}\}$ is a basis of $P_{m-1}$.

This linear system is solvable because $h(x)$ is conditionally positive definite(c.p.d.) of order $m=max\{ \lceil \frac{\beta}{2}\rceil , 0\}$ where $\lceil \frac{\beta}{2}\rceil \}$ denotes the smallest integer greater than or equal to $\frac{\beta}{2}$. 

Besides the linear system, another important object is the function space. Each function of the form (1) induces a function space called {\bf native space} denoted by ${\cal C}_{h,m}(R^{n})$, abbreviated as ${\cal C}_{h,m}$, where $m$ denotes its order of conditional positive definiteness. For each member $f$ of ${\cal C}_{h,m}$ there is a seminorm $\|f\|_{h}$, called the $h$-norm of $f$. The definition and characterization of the native space can be found in \cite{Lu1}, \cite{Lu2}, \cite{Lu3-1}, \cite{MN1}, \cite{MN2} and \cite{We}. In this paper all interpolated functions belong to the native space.

Although our interpolated functions are defined in the entire $R^{n}$, interpolation will occur in a simplex. The definition of simplex can be found in \cite{Fl}. A 1-simplex is a line segment, a 2-simplex is a triangle, and a 3-simplex is a tetrahedron with four vertices.

Let $T_{n}$ be an n-simplex in $R^{n}$ and $v_{i},\ 1\leq i\leq n+1$ be its vertices. Then any point $x\in T_{n}$ can be written as convex combination of the vertices:
$$x=\sum_{i=1}^{n+1}c_{i}v_{i},\ \sum_{i=1}^{n+1}c_{i}=1,\ c_{i}\geq 0.$$
The numbers $c_{1},\cdots ,c_{n+1}$ are called the barycentric coordinates of $x$. For any n-simplex $T_{n}$, the {\bf evenly spaced points} of degree $l$ are those points whose barycentric coordinates are of the form
$$(\frac{k_{1}}{l},\frac{k_{2}}{l},\cdots,\frac{k_{n+1}}{l}),\ k_{i}\ nonnegative\ integers\ with\ \sum_{i=1}^{n+1}k_{i}=l.$$
It's easily seen that the number of evenly spaced points of degree $l$ in $T_{n}$ is exactly $$N=dimP_{l}^{n}=\left( \begin{array}{c}
                        n+l \\
                        n

                      \end{array} \right) $$
where $P_{l}^{n}$ denotes the space of polynomials of degree not exceeding $l$ in n variables. Moreover, such points form a determining set for $P_{l}^{n}$, as is shown in \cite{Bo}.

In this paper the evaluation argument $x$ will be a point in an n-simplex, and the set $X$ of centers will be the evenly spaced points in that n-simplex.

\section{Fundamental Theory}
Before introducing the main theorem, we need to define two constants.
\begin{de}
Let $n$ and $\beta$ be as in (1). The numbers $\rho$ and $\Delta_{0}$ are defined as follows.
\begin{list}
  {(\alph{bean})}{\usecounter{bean} \setlength{\rightmargin}{\leftmargin}}
  \item Suppose $\beta <n-3$. Let $s=\lceil \frac{n-\beta -3}{2}\rceil $. Then 
    \begin{list}{(\roman{milk})}{\usecounter{milk} \setlength{\rightmargin}{\leftmargin}}
      \item if $\beta <0,\ \rho=\frac{3+s}{3}\ and\  \Delta_{0}=\frac{(2+s)(1+s)\cdots 3}{
 \rho^{2}};$
      \item if $\beta >0,\ \rho=1+\frac{s}{2\lceil \frac{\beta}{2}\rceil +3} \ and \ \Delta_{0}=\frac{(2m+2+s)(2m+1+s)\cdots (2m+3)}{\rho^{2m+2}}$ \\
where $ m=\lceil \frac{\beta}{2}\rceil$.          
    \end{list}
  \item Suppose $n-3\leq \beta <n-1$. Then $\rho=1$ and $\Delta_{0}=1$.
  \item Suppose $\beta \geq n-1$. Let $s=-\lceil \frac{n-\beta -3}{2}\rceil $. Then
 $$\rho =1\ and \ \Delta_{0}=\frac{1}{(2m+2)(2m+1)\cdots (2m-s+3)} \ where \ m=\lceil \frac{\beta}{2}\rceil.$$  
\end{list}
\end{de}

The following theorem is the cornerstone of our theory. We cite it directly from \cite{Lu3} with a slight modification to make it easier to understand.

\begin{thm}
  Let $h$ be as in (1). For any positive number $b_{0}$, let $C=\max \left\{ \frac{2}{3b_{0}},8\rho\right\}$ and $\delta_{0}=\frac{1}{3C}$. For any n-simplex $Q$ of diameter $r$ satisfying $\frac{1}{3C}\leq r\leq \frac{2}{3C}$(note that $\frac{2}{3C}\leq b_{0}$), if $f\in {\cal C}_{h,m}$,
\begin{eqnarray}
  |f(x)-s(x)|\leq 2^{\frac{n+\beta-7}{4}}\pi^{\frac{n-1}{4}}\sqrt{n\alpha_{n}}c^{\frac{\beta}{2}-l}\sqrt{\Delta_{0}}\sqrt{3C}\sqrt{\delta}(\lambda')^{\frac{1}{\delta}}\|f\|_{h}
\end{eqnarray}
holds for all $x\in Q$ and $0<\delta<\delta_{0}$, where $s(x)$ is defined as in (2) with $x_{1},\cdots ,x_{N}$ the evenly spaced points of degree $l$ in $Q$ satisfying $\frac{1}{3C\delta}\leq l\leq \frac{2}{3C\delta}$. The constant $\alpha_{n}$ denotes the volume of the unit ball in $R^{n}$, and $0<\lambda'<1$ is given by 
$$\lambda'=\left(\frac{2}{3}\right)^{\frac{1}{3C}}$$
which only in some cases mildly depends on the dimension n. 
\end{thm}
{\bf Remark}:(a)Note that the right-hand side of (4) approaches zero as $\delta\rightarrow 0^{+}$. This is the key to understanding Theorem2.2. The number $\delta$ is in spirit equivalent to the well-known fill-distance. Although the centers $x_{1},\cdots,x_{N}$ are not purely scattered, the shape of the simplex is controlled by us. Hence the distribution of the centers is practically quite flexible. (b)In (4) the shape parameter c plays a crucial role and greatly influences the error bound. This provides us with a theoretical ground of choosing the optimal c. However we need further work before presenting useful criteria.

In this paper all interpolated functions belong to a kind of space defined as follows.
\begin{de}
For any positive number $\sigma$,
$$E_{\sigma}:=\left\{ f\in L^{2}(R^{n}):\ \int |\hat{f}(\xi)|^{2}e^{\frac{|\xi|^{2}}{\sigma}}d\xi<\infty \right\}$$
where $\hat{f}$ denotes the Fourier transform of $f$. For each $f\in E_{\sigma}$, its norm is 
$$\|f\|_{E_{\sigma}}:=\left\{ \int|\hat{f}(\xi)|^{2}e^{\frac{|\xi|^{2}}{\sigma}}d\xi\right\}^{1/2}$$.
\end{de}
The following lemma is cited from \cite{Lu5}.
\begin{lem}
Let $h$ be as in (1). For any $\sigma>0$, if $\beta<0$, $|n+\beta|\geq 1$ and $n+\beta+1\geq 0$, then $E_{\sigma}\subseteq {\cal C}_{h,m}(R^{n})$ and for any $f\in E_{\sigma}$, the seminorm $\|f\|_{h}$ of $f$ satisfies 
$$\|f\|_{h}\leq 2^{-n-\frac{1+\beta}{4}}\pi^{-n-\frac{1}{4}}c^{\frac{1-n-\beta}{4}}\left\{ (\xi^{*})^{\frac{n+\beta+1}{2}}e^{c\xi^{*}-\frac{(\xi^{*})^{2}}{\sigma}}\right\}^{1/2}\|f\|_{E_{\sigma}}$$
where $$\xi^{*}:=\frac{c\sigma+\sqrt{c^{2}\sigma^{2}+4\sigma(n+\beta+1)}}{4}$$.
\end{lem}
\begin{cor}
Under the conditions of Theorem2.2, if $f\in E_{\sigma},\ \beta<0,\ |n+\beta|\geq 1$ and $n+\beta+1\geq 0$, (4) can be transformed into
\begin{eqnarray}
|f(x)-s(x)|\leq 2^{-\frac{3n}{4}-2}\pi^{-\frac{3}{4}n-\frac{1}{2}}\sqrt{n\alpha_{n}}\sqrt{\Delta_{0}}\sqrt{3C}c^{\frac{\beta-n+1-4l}{4}}\left\{(\xi^{*})^{\frac{n+\beta+1}{2}}e^{c\xi^{*}-\frac{(\xi^{*})^{2}}{\sigma}}\right\}^{1/2}\sqrt{\delta}(\lambda')^{\frac{1}{\delta}}\|f\|_{E_{\sigma}}
\end{eqnarray}
where $$\xi^{*}:=\frac{c\sigma+\sqrt{c^{2}\sigma^{2}+4\sigma(n+\beta+1)}}{4}$$.
\end{cor}
{\bf Proof}. This is an immediate result of Theorem2.2 and Lemma2.4. \hspace{5cm} $\sharp$\\
\\
Note that Corollary2.5 covers the very useful case $\beta=-1,\ n\geq 2$. However the case $\beta=-1,\ n=1$ is excluded. For this case we need a different approach.
\begin{lem}
Let $\sigma>0,\ \beta=-1$ and $n=1$. For any $f\in E_{\sigma}$,
$$\|f\|_{h}\leq 2^{-(n+\frac{1}{4})}\pi^{-1}\left\{ \frac{1}{ln2}+2\sqrt{3}M(c)\right\}^{1/2}\|f\|_{E_{\sigma}}$$
where $M(c):=e^{1-\frac{1}{c^{2}\sigma}}$ if $c\leq \frac{2}{\sqrt{3\sigma}}$ and $M(c):=g(\frac{c\sigma+\sqrt{c^{2}\sigma^{2}+4\sigma}}{4})$ if $c>\frac{2}{\sqrt{3\sigma}}$, where $g(\xi):=\sqrt{c\xi}e^{c\xi-\frac{\xi^{2}}{\sigma}}$.
\end{lem}
{\bf Proof}. This is just Theorem2.5 of \cite{Lu5}. \hspace{9cm} $\sharp$

\begin{cor}
Let $\sigma>0,\ \beta=-1$ and $n=1$. Under the conditions of Theorem2.2, if $f\in E_{\sigma}$, (4) can be transformed into 
\begin{eqnarray}
|f(x)-s(x)|\leq 2^{\frac{\beta-3n}{4}-2}\pi^{\frac{n-5}{4}}\sqrt{n\alpha_{n}}\sqrt{\Delta_{0}}\sqrt{3C}c^{\frac{\beta}{2}-l}\left\{ \frac{1}{ln2}+2\sqrt{3}M(c)\right\}^{1/2}\sqrt{\delta}(\lambda')^{\frac{1}{\delta}}\|f\|_{E_{\sigma}}
\end{eqnarray}
where $M(c)$ is defined as in Lemma2.6.
\end{cor}
{\bf Proof}. This is an immediate result of Theorem2.2 and Lemma2.6. \hspace{4.8cm} $\sharp$\\
\\
Now we have dealt with the most useful cases for $\beta<0$. The next step is to treat $\beta>0$.
\begin{lem}
 Let $\sigma>0,\ \beta>0$ and $n\geq 1$. For any $f\in E_{\sigma}$,
$$\|f\|_{h}\leq d_{0}c^{\frac{1-\beta-n}{4}}\left\{ \frac{(\xi^{*})^{\frac{1+\beta+n}{2}}e^{c\xi^{*}}}{e^{\frac{(\xi^{*})^{2}}{\sigma}}}\right\} ^{1/2}\|f\|_{E_{\sigma}}$$
where $\xi^{*}=\frac{c\sigma+\sqrt{c^{2}\sigma^{2}+4\sigma(1+\beta+n)}}{4}$ and $d_{0}$ is a constant depending on $n,\ \beta$ only.
\end{lem}
{\bf Proof}. This is just Theorem2.8 of \cite{Lu5}. \hspace{8.9cm} $\sharp$

\begin{cor}
  Let $\sigma>0,\ \beta>0$ and $n\geq 1$. If $f\in E_{\sigma}$, (4) can be transformed into 
\begin{eqnarray}
|f(x)-s(x)|\leq 2^{\frac{n+\beta-7}{4}}\pi^{\frac{n-1}{4}}\sqrt{n\alpha_{n}}\sqrt{\Delta_{0}}\sqrt{3C}d_{0}c^{\frac{1+\beta-n-4l}{4}}\left\{ \frac{(\xi^{*})^{\frac{1+\beta+n}{2}}e^{c\xi^{*}}}{e^{\frac{(\xi^{*})^{2}}{\sigma}}}\right\} ^{1/2}\sqrt{\delta}(\lambda')^{\frac{1}{\delta}}\|f\|_{E_{\sigma}}
\end{eqnarray}
where $d_{0},\ \xi^{*}$ are as in Lemma2.8.
\end{cor}
{\bf Proof}. This is an immediate result of Theorem2.2 and Lemma2.8. \hspace{4.5cm} $\sharp$

\section{Criteria of Choosing c}
Note that in (5),(6) and (7), there is a main function of c. As in \cite{Lu5}, let's call this function the MN function, denoted by $MN(c)$, and its graph the MN curve. The optimal choice of c is then the number minimizing $MN(c)$. However, unlike \cite{Lu5}, the range of c is the entire interval $(0,\infty)$, rather than a proper subset of $(0,\infty)$.

We now begin our criteria.\\
\\
{\bf Case1}. \fbox{$\beta<0,\ |n+\beta|\geq 1$ and $n+\beta+1\geq 0$} Let $f\in E_{\sigma}$ and $h$ be as in (1). Under the conditions of Theorem2.2, for any fixed $\delta$ satisfying $0<\delta<\delta_{0}$, the optimal value of c in $(0,\infty)$ is the number minimizing
$$MN(c):=c^{\frac{\beta-n+1-4l}{4}}\left\{(\xi^{*})^{\frac{n+\beta+1}{2}}e^{c\xi^{*}-\frac{(\xi^{*})^{2}}{\sigma}}\right\}^{1/2}$$
where $$\xi^{*}=\frac{c\sigma+\sqrt{c^{2}\sigma^{2}+4\sigma(n+\beta+1)}}{4}$$.\\
\\
{\bf Reason}: This is a direct consequence of (5). \hspace{8cm} $\sharp$\\
\\
{\bf Remark}:(a)It's easily seen that $MN(c)\rightarrow\infty$ as $c\rightarrow\infty$. Also, if $n+\beta+1>0,\ MN(c)\rightarrow\infty$ as $c\rightarrow 0^{+}$. (b)Case1 covers the frequently seen case $\beta=-1,\ n\geq 2$. (c)The number c minimizing $MN(c)$ can be easily found by Mathematica or Matlab.\\
\\
{\bf Numerical Results}:\\
\\
\begin{figure}[h]
\centering
\includegraphics[scale=1.0]{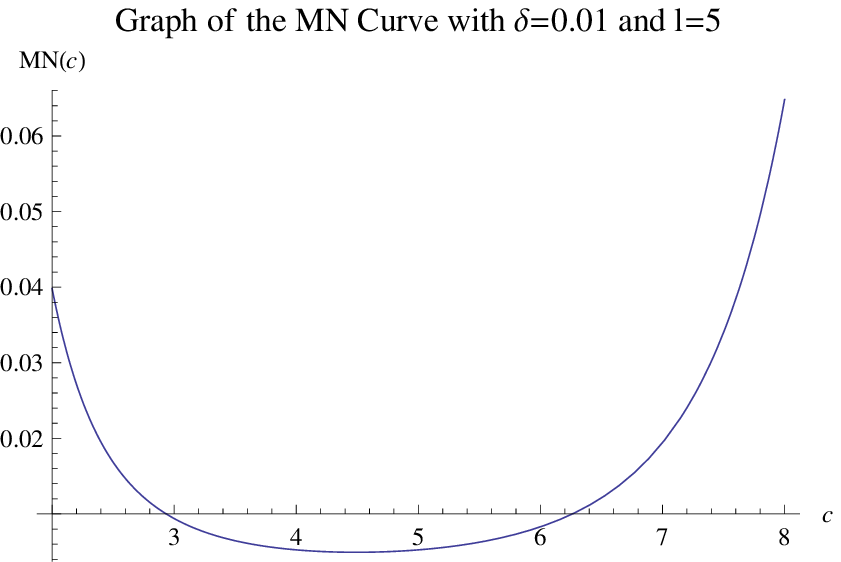}
\caption{Here $n=2,\beta=-1,\sigma=1$ and $b_{0}=1$.}
\end{figure}

\clearpage

\begin{figure}[t]
\centering
\includegraphics[scale=1.0]{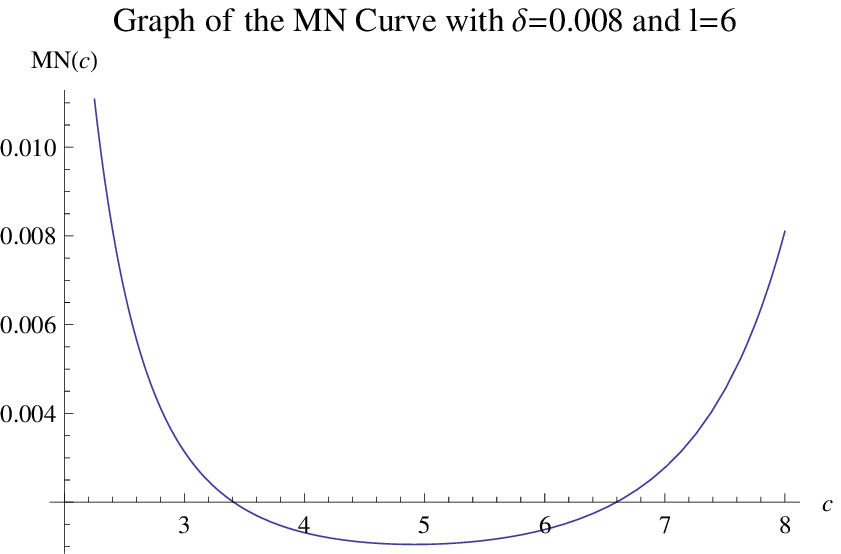}
\caption{Here $n=2,\beta=-1,\sigma=1$ and $b_{0}=1$.}

\includegraphics[scale=1.0]{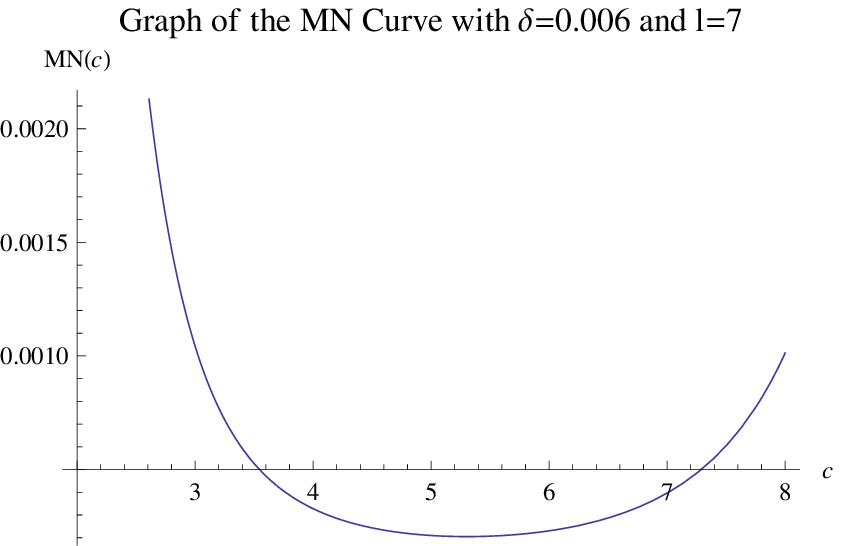}
\caption{Here $n=2,\beta=-1,\sigma=1$ and $b_{0}=1$.}

\includegraphics[scale=1.0]{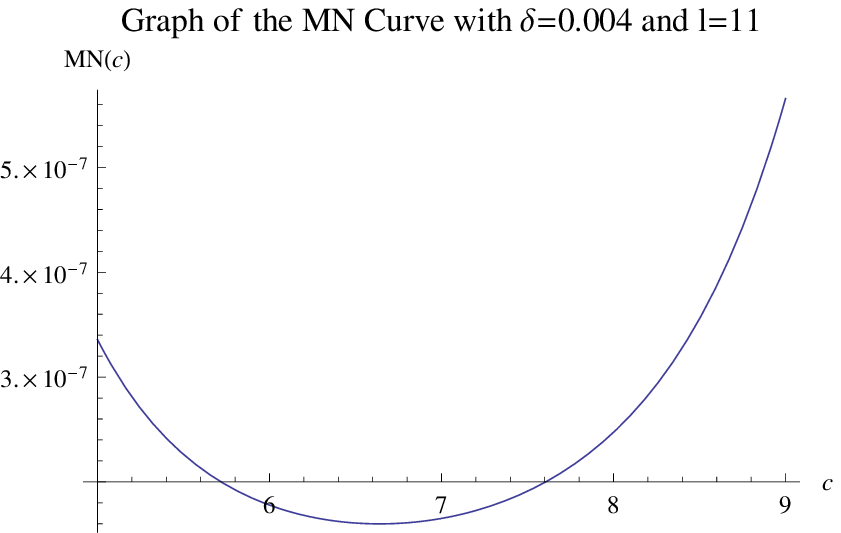}
\caption{Here $n=2,\beta=-1,\sigma=1$ and $b_{0}=1$.}

\end{figure}

\clearpage

\begin{figure}[t]
\centering
\includegraphics[scale=1.0]{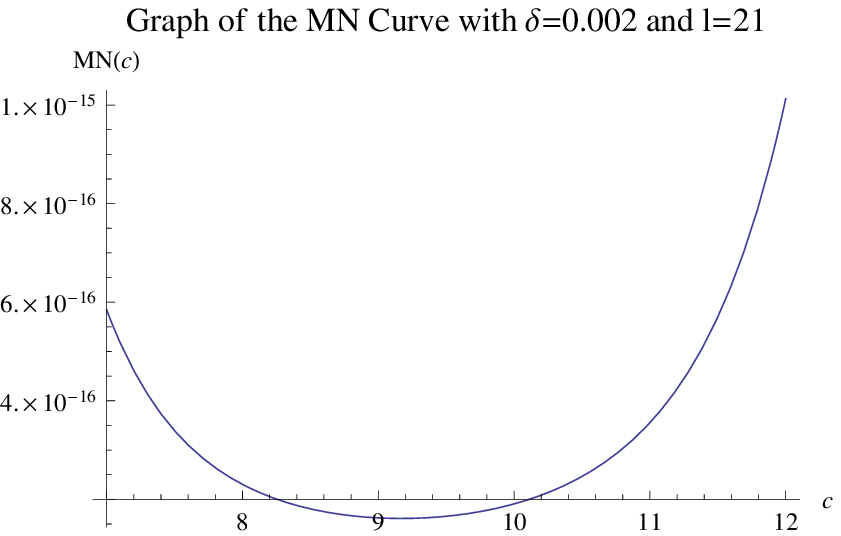}
\caption{Here $n=2,\beta=-1,\sigma=1$ and $b_{0}=1$.}
\end{figure}

{\bf Case2}. \fbox{$\beta=-1$ and $n=1$} Let $f\in E_{\sigma}$ and $h$ be as in (1). Under the conditions of Theorem2.2, for any fixed $\delta$ satisfying $0<\delta<\delta_{0}$, the optimal value of c in $(0,\infty)$ is the number minimizing
$$MN(c):=c^{\frac{\beta}{2}-l}\left\{\frac{1}{ln2}+2\sqrt{3}M(c)\right\}^{1/2}$$
where
$$M(c):=\left\{ \begin{array}{ll}
                 e^{1-\frac{1}{c^{2}\sigma}}     & \mbox{if $0<c\leq \frac{2}{\sqrt{3\sigma}}$,} \\
                 g(\frac{c\sigma+\sqrt{c^{2}\sigma^{2}+4\sigma}}{4})                                                & \mbox{if $\frac{2}{\sqrt{3\sigma}}<c$}

                \end{array} \right. $$, $g$ being defined by $g(\xi):=\sqrt{c\xi}e^{c\xi-\frac{\xi^{2}}{\sigma}}$.\\
\\
{\bf Reason}: This is a direct result of (6). \hspace{9cm} $\sharp$\\
\\
{\bf Remark}: Note that $MN(c)\rightarrow \infty$ both as $c\rightarrow \infty$ and $c\rightarrow 0^{+}$. Now let's see some numerical examples.

\begin{figure}[h]
\centering
\includegraphics[scale=0.8]{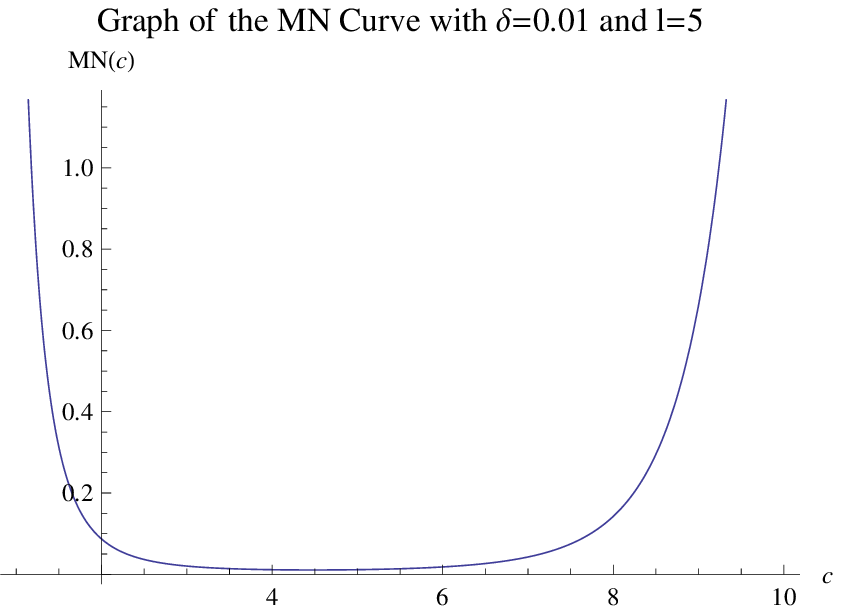}
\caption{Here $n=1,\beta=-1,\sigma=1$ and $b_{0}=1$.}

\end{figure}

\clearpage

\begin{figure}[t]
\centering
\includegraphics[scale=1.0]{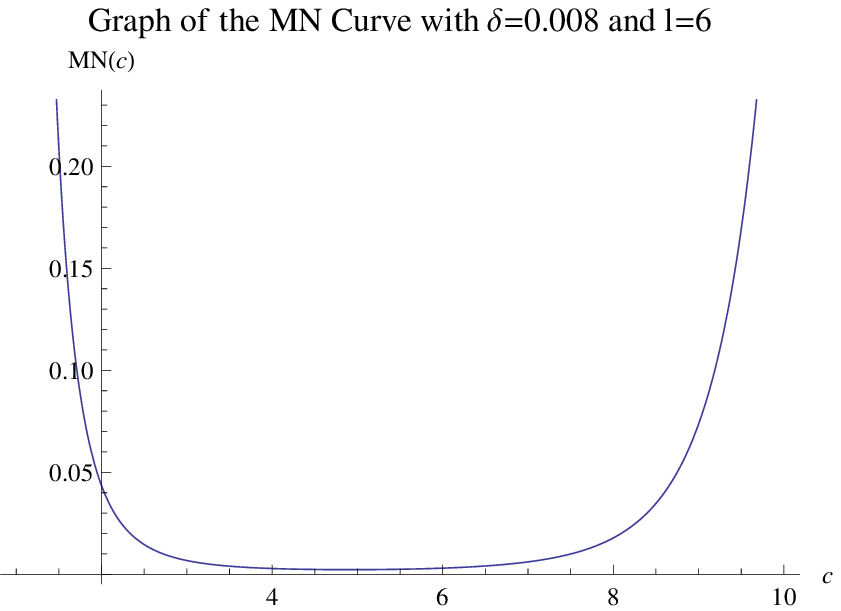}
\caption{Here $n=1,\beta=-1,\sigma=1$ and $b_{0}=1$.}

\includegraphics[scale=1.0]{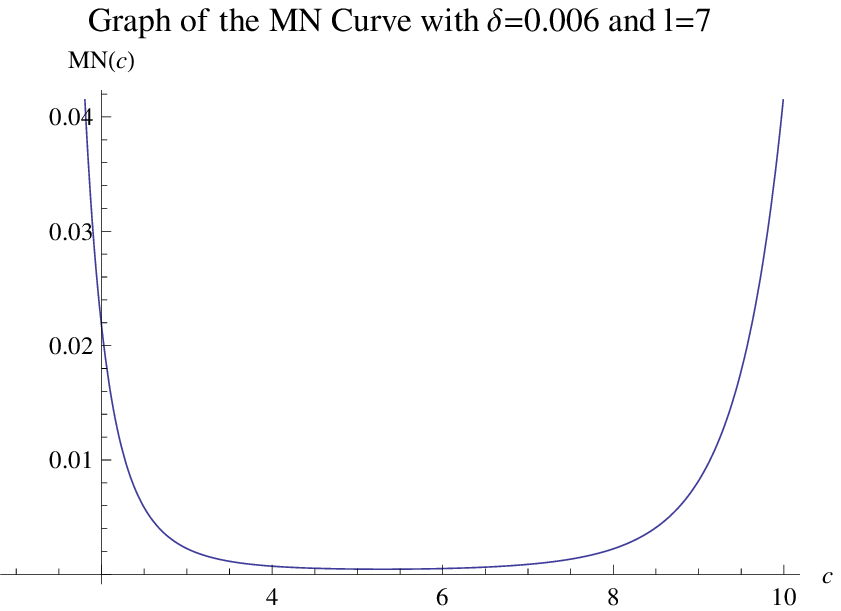}
\caption{Here $n=1,\beta=-1,\sigma=1$ and $b_{0}=1$.}

\includegraphics[scale=1.0]{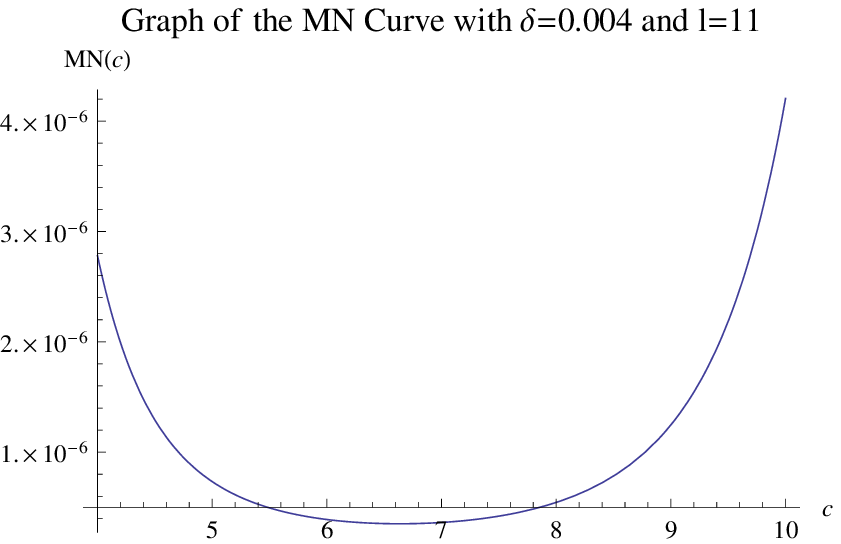}
\caption{Here $n=1,\beta=-1,\sigma=1$ and $b_{0}=1$.}

\end{figure}
\clearpage
\begin{figure}[t]
\centering
\includegraphics[scale=1.0]{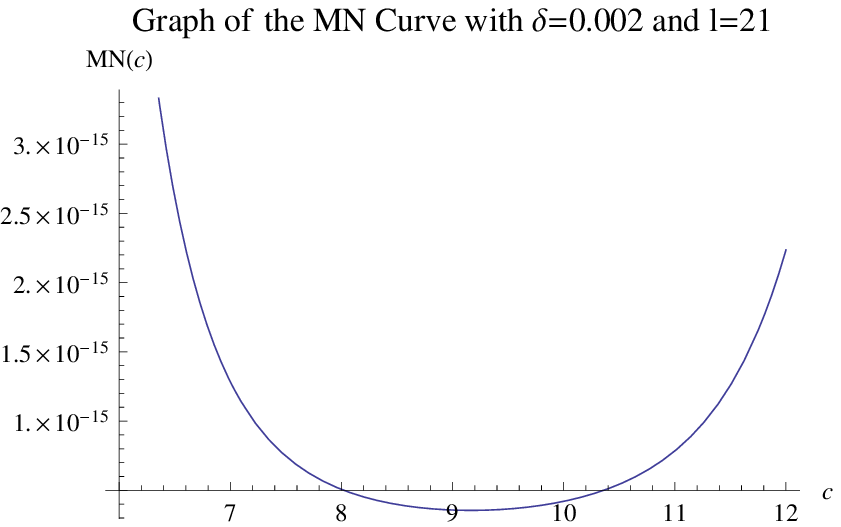}
\caption{Here $n=1,\beta=-1,\sigma=1$ and $b_{0}=1$.}
\end{figure}

{\bf Case3}. \fbox{$\beta>0$ and $n\geq 1$} Let $f\in E_{\sigma}$ and $h$ be as in (1). Under the conditions of Theorem2.2, for any fixed $\delta$ satisfying $0<\delta<\delta_{0}$, the optimal value of c in $(0,\infty)$ is the number minimizing
$$MN(c):=c^{\frac{1+\beta-n-4l}{4}}\left\{\frac{(\xi^{*})^{\frac{1+\beta+n}{2}}e^{c\xi^{*}}}{e^{\frac{(\xi^{*})^{2}}{\sigma}}}\right\}^{1/2}$$
, where
$$\xi^{*}=\frac{c\sigma+\sqrt{c^{2}\sigma^{2}+4\sigma(1+\beta+n)}}{4}$$.\\
\\
{\bf Reason}: This follows from (7). \hspace{10cm} $\sharp$\\
\\
{\bf Remark}: By observing that
$$c\xi^{*}-\frac{(\xi^{*})^{2}}{\sigma}=\frac{1}{16}\left[ 2c^{2}\sigma+2c\sqrt{c^{2}\sigma^{2}+4\sigma(n+\beta+1)}-(4n+\beta+1)\right]$$
, we can easily obtain useful results as follows. (a)If $1+\beta-n-4l>0$, $\lim_{c\rightarrow0^{+}}MN(c)=0$. (b)If $1+\beta-n-4l<0,\ \lim_{c\rightarrow 0^{+}}MN(c)=\infty$. (c)If $1+\beta-n-4l=0,\ \lim_{c\rightarrow 0^{+}}MN(c)$ is a finite positive number. (d)$\lim_{c\rightarrow \infty}MN(c)=\infty$.\\
\\
{\bf Numerical Results}: For simplicity, we offer results for $n=1$ only. In fact for $n\geq 1$ similar results can be presented without slight difficulty.

\begin{figure}[t]
\centering
\includegraphics[scale=1.0]{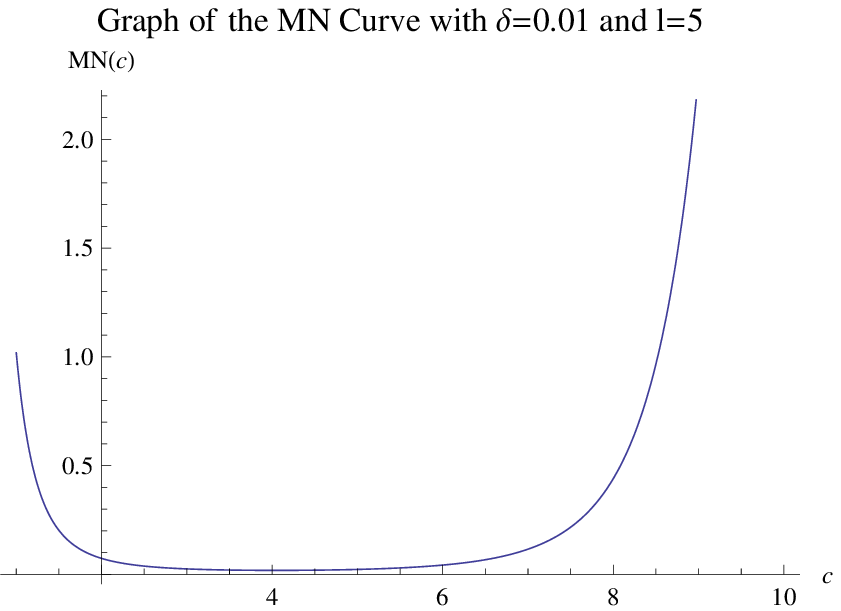}
\caption{Here $n=1,\beta=1,\sigma=1$ and $b_{0}=1$.}

\includegraphics[scale=1.0]{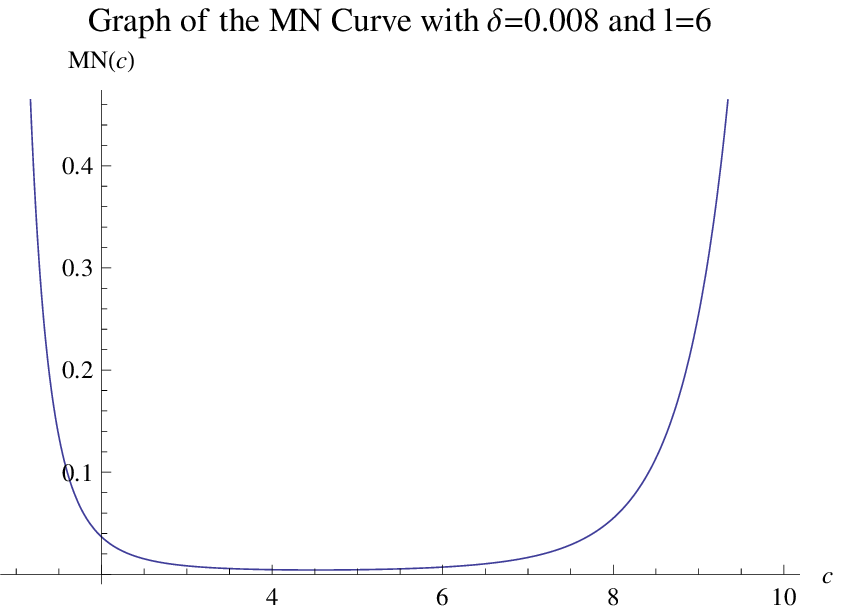}
\caption{Here $n=1,\beta=1,\sigma=1$ and $b_{0}=1$.}

\includegraphics[scale=1.0]{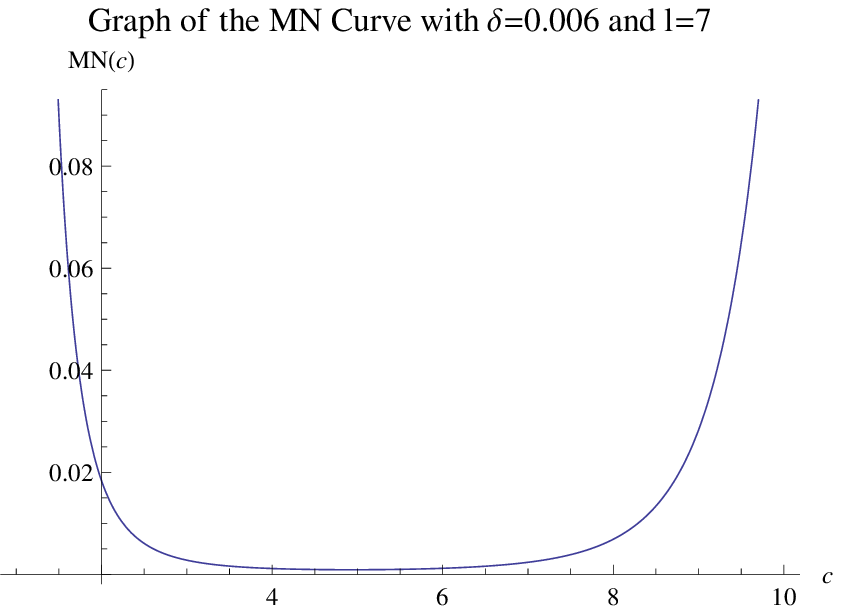}
\caption{Here $n=1,\beta=1,\sigma=1$ and $b_{0}=1$.}

\end{figure}
\clearpage

\begin{figure}[t]
\centering
\includegraphics[scale=1.0]{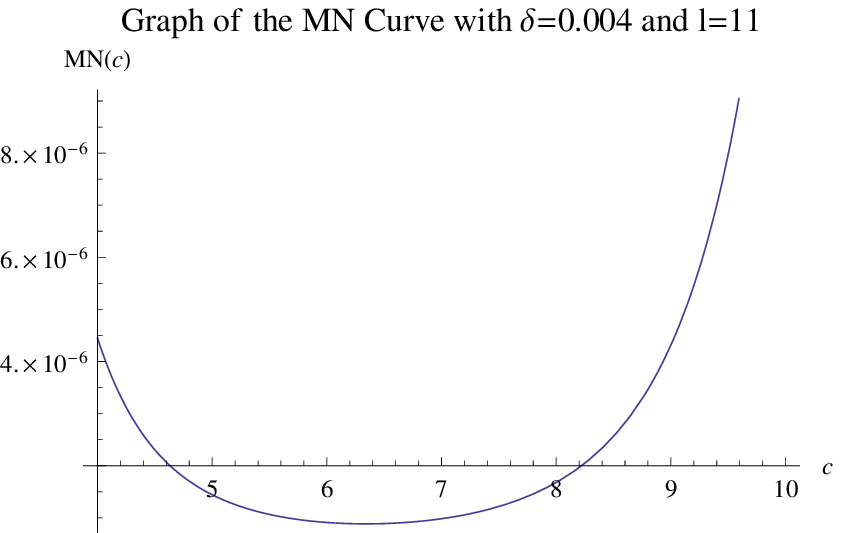}
\caption{Here $n=1,\beta=1,\sigma=1$ and $b_{0}=1$.}

\includegraphics[scale=1.0]{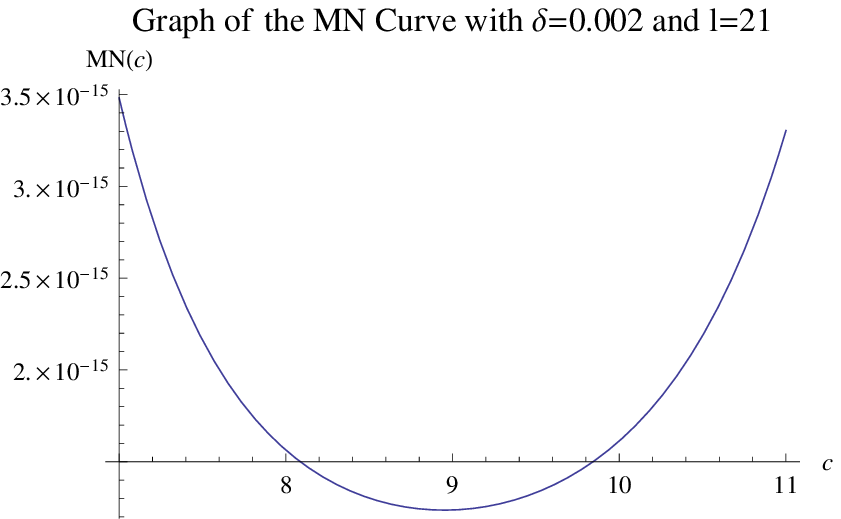}
\caption{Here $n=1,\beta=1,\sigma=1$ and $b_{0}=1$.}

\end{figure}



\begin{thebibliography}{99}
\bibitem{Ab}Abramowitz and Segun,
{\em Handbook of Mathematical Functions,}
Dover Publications, INC., New York. 

\bibitem{Bo}L.P. Bos,
{\em Bounding the Lebesgue function for Lagrange interpolation in a simplex,}
J. Approx. Theory, 38(1983)43-59.

\bibitem{Fl}W. Fleming,
{\em Functions of Several Variables, Second Edition,}
Springer-Verlag, 1977.

\bibitem{Lu1}L.T. Luh,
{\em The Equivalence Theory of Native Spaces,}
Approx. Theory Appl. (2001), 17:1, 76-96.

\bibitem{Lu2}L.T. Luh,
{\em The Embedding Theory of Native Spaces,}
Approx. Theory Appl. (2001), 17:4, 90-104.

\bibitem{Lu3}L.T. Luh,
{\em An Improved Error Bound for Multiquadric Interpolation,}
Inter. J. Numeric. Methods Appl. Vol. 1, No.2, pp. 101-120, 2009.

\bibitem{Lu3-1}L.T. Luh,
{\em On Wu and Schaback's Error Bound,}
Inter. J. Numeric. Methods Appl. Vol. 1, No2, pp. 155-174, 2009. 

\bibitem{Lu4}L.T. Luh,
{\em The Mystery of the Shape Parameter,}
Math ArXiv.

\bibitem{Lu5}L.T. Luh,
{\em The Mystery of the Shape Parameter II,}
Math ArXiv.

\bibitem{Lu6}L.T. Luh,
{\em The Mystery of the Shape Parameter III,}
Math ArXiv.

\bibitem{MN1}W.R. Madych and S.A. Nelson,
{\em Multivariate interpolation and conditionally positive definite function,}
Approx. Theory Appl. 4, No. 4(1988), 77-89.

\bibitem{MN2}W.R. Madych and S.A. Nelson,
{\em Multivariate interpolation and conditionally positive definite function, II,}
Math. Comp. 54(1990), 211-230.

\bibitem{MN3}W.R. Madych,
{\em Miscellaneous Error Bounds for Multiquadric and Related Interpolators,}
Computers Math. Applic. Vol. 24, No. 12, pp. 121-138, 1992.
 
\bibitem{We}H. Wendland,
{\em Scattered Data Approximation,}
Cambridge University Press, (2005).


\end{thebibliography}
\end{document}